# Exact solution of a nonlinear optimization problem

**Zhixi Bian[1], Hongyun Wang[2]\*, Qiaoer Zhou[1], and Ali Shakouri[1]**

[1] Department of Electrical Engineering, Jack Baskin School of Engineering, University of California, Santa Cruz, CA 95064, USA

[2] Department of Applied Mathematics and Statistics, Jack Baskin School of Engineering, University of California, Santa Cruz, CA 95064, USA

\* Corresponding author, email: hongwang@soe.ucsc.edu

## Abstract

In this manuscript, we solve a nonlinear optimization problem in the study of maximizing cooling temperature using inhomogeneous thermoelectric materials.

## 1  The optimization problem

We consider the idealized case where the Seebeck coefficient and the electrical conductivity can be varied while the value of $ZT^2$ remains constant. For such an ideal graded material, the maximum cooling temperature (with respect to the electrical current) is given by

$$\Delta T_{\max} = \frac{1}{2} ZT^2 \frac{\int_0^L S(x)dx \int_0^x S(x')dx'}{\int_0^L dx \int_0^x S^2(x')dx'}$$

where $ZT^2$ is a constant, independent of temperature, $S(x)$ is the profile of Seebeck coefficient, $x$ is the one-dimensional coordinate along the material and L is the length of the material. Using the normalized coordinate $\tilde{x} = \frac{x}{L}$, letting $\tilde{S}(\tilde{x}) = S(x)$, and rewriting the two double integrals, we obtain

$$\Delta T_{max} = \frac{1}{2} Z T^2 F[\tilde{S}(\tilde{x})], \qquad F[\tilde{S}(\tilde{x})] \equiv \frac{\frac{1}{2}\left(\int_0^1 S(\tilde{x}) d\tilde{x}\right)^2}{\int_0^1 (1-\tilde{x})\tilde{S}^2(\tilde{x}) d\tilde{x}}$$

For the mathematical convenience, we drop the tilde and go back to the notations of $x$ and $S(x)$. The mathematical problem is to maximize the functional

$$F[S(x)] \equiv \frac{\frac{1}{2}\left(\int_0^1 S(x) dx\right)^2}{\int_0^1 (1-x) S^2(x) dx}$$

where $S(x)$ is the unknown profile of Seebeck coefficient that is subject to the constraint

$$S_0 \leq S(x) \leq S_1, \quad x \in [0, 1]$$

## 2  Solution of the optimization problem

Below we will show that the optimal Seebeck profile is given by

$$S_{opt}(x) = \begin{cases} S_0, & 0 \leq x \leq \frac{1}{2} \\ \frac{S_0}{2(1-x)}, & \frac{1}{2} < x < 1 - \frac{S_0}{2S_1} \\ S_1, & 1 - \frac{S_0}{2S_1} \leq x \leq 1 \end{cases}$$

It is straightforward to verify that the maximum cooling temperature for $S_{opt}(x)$ is

$$\Delta T_{max} = \left(\frac{1}{2} Z T^2\right)\left(1 + \frac{1}{2}\ln\left(\frac{S_1}{S_0}\right)\right)$$

Now we show step by step the mathematical analysis that leads to $S_{opt}(x)$.

**Step #1: The optimal S(x) must be non-decreasing.**

Suppose we discretize $S(x)$ as

$$S(x) = y_j \quad \text{for } x \in \left(\frac{j}{N}, \frac{j+1}{N}\right)$$

The numerator and denominator of $F[S(x)]$ are expressed as

$$\frac{1}{2}\left(\int_0^1 S(x)dx\right)^2 = \frac{1}{2}\left(\frac{1}{N}\sum_{j=0}^{N-1} y_j\right)^2$$

$$\int_0^1 S^2(x)(1-x)dx = \frac{1}{N}\sum_{j=0}^{N-1}\left(1-\frac{j}{N}\right)y_j^2$$

If $S(x)$ is not non-decreasing, then we can find a pair of indices $(j, k)$ such that $j < k$ but $y_j > y_k$, and by exchanging the values of $y_j$ and $y_k$, we have

$$\frac{1}{2}\left(\frac{1}{N}\sum_{j=0}^{N-1} y_j\right)^2 \quad \text{is unchanged}$$

but $\frac{1}{N}\sum_{j=0}^{N-1}\left(1-\frac{j}{N}\right)y_j^2$ is decreased.

It follows that the value of $F[S(x)]$ is increased after the exchange. Therefore, the optimal $S(x)$ must be non-decreasing.

**Step #2: The optimal S(x) has 3 segments**

$$S_{opt}(x) = \begin{cases} S_0, & 0 \leq x \leq x_0 \\ S_0 < S(x) < S_1, & x_0 < x < x_1 \\ S_1, & x_1 \leq x \leq 1 \end{cases}$$

This follows directly from that $S_{opt}(x)$ must be non-decreasing (result of step #1).

**Step #3: The middle segment of $S_{opt}(x)$ must satisfy**

$$S_{opt}(x) = \frac{q}{1-x} \quad \text{for } x \in (x_0, x_1)$$

That is, $S_{opt}(x)$ must have the shape as shown in Figure A1 below.

Let $\Delta S(x)$ be a function that is non-zero only in the middle segment. That is, $\Delta S(x) = 0$ in $[0, x_0]$ and in $[x_1, 1]$. Consider a small perturbation to the middle segment of $S_{opt}(x)$:

$$S_{opt}(x) + \varepsilon \Delta S(x)$$

When $\varepsilon$ is small enough, $S_{opt}(x) + \varepsilon \Delta S(x)$ is between $S_0$ and $S_1$ (i.e. satisfying the constraint of the optimization).

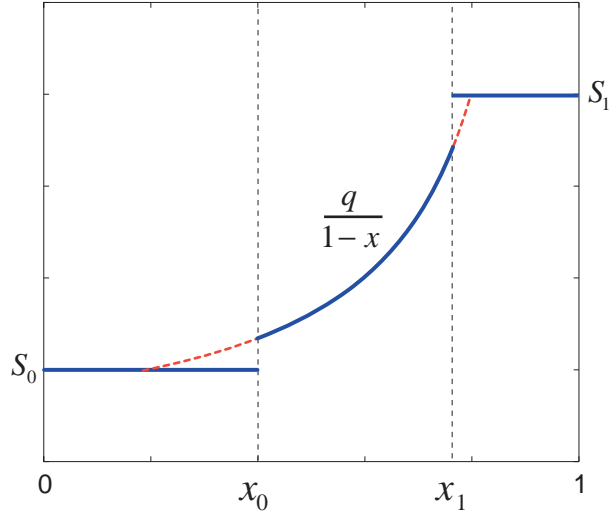

Figure A1

$S_{opt}(x)$ is the optimal Seebeck profile implies

$$F[S_{opt}(x) + \varepsilon \Delta S(x)] \le F[S_{opt}(x)] \qquad \text{for } \varepsilon \text{ small enough}$$

which leads to

$$\left. \frac{dF[S_{opt}(x) + \varepsilon \Delta S(x)]}{d\varepsilon} \right|_{\varepsilon=0} = 0$$

Restricting our attention to perturbations satisfying $\int_{x_0}^{x_1} \Delta S(x) dx = 0$, we have

$$\left. \frac{dF[S_{opt}(x)+\varepsilon \Delta S(x)]}{d\varepsilon} \right|_{\varepsilon=0} = \frac{-2\left(\int_0^1 S_{opt}(x)dx\right)^2}{\left(\int_0^1 (1-x)S_{opt}^2(x)dx\right)^2} \int_0^1 (1-x)S_{opt}(x)\Delta S(x)dx$$

$$\Longrightarrow \quad \int_{x_0}^{x_1} (1-x)S_{opt}(x)\Delta S(x)dx = 0 \qquad \text{for all } \Delta S(x) \text{ satisfying } \int_{x_0}^{x_1} \Delta S(x)dx = 0$$

$$\Longrightarrow \quad (1-x)S_{opt}(x) = \text{const} \quad \text{for } x \in (x_0, x_1)$$

$$\Longrightarrow \quad S_{opt}(x) = \frac{q}{1-x} \quad \text{for } x \in (x_0, x_1)$$

**Step #4: $S_{opt}(x)$ must be continuous at both $x_0$ and $x_1$**

We show the continuity using the method of <u>proof by contradiction.</u>

Suppose $S_{opt}(x)$ is discontinuous at $x_0$. Since $S_{opt}(x)$ is non-decreasing, we have

$$S_0 < \frac{q}{1-x_0} \quad \text{and} \quad \frac{q}{1-x_1} \leq S_1$$

$$\Longrightarrow \quad (1-x_0)S_0 < q \leq (1-x_1)S_1$$

$$\Longrightarrow \quad (1-x_0)S_0 - (1-x_1)S_1 < 0 \qquad\qquad\qquad (\text{A01})$$

Let $\Delta S(x)$ be a function that is non-zero only in $[x_0 - \delta, x_0]$ and $[x_1, x_1 + \delta]$. Specifically

$$\Delta S(x) = \begin{cases} 1, & x \in [x_0 - \delta, x_0] \\ -1, & x \in [x_1, x_1 + \delta] \\ 0, & \text{otherwise} \end{cases}$$

Consider a small perturbation to $S_{opt}(x)$:

$$S_{opt}(x) + \varepsilon \Delta S(x)$$

When $\varepsilon$ is <u>positive and small enough</u>, $S_{opt}(x) + \varepsilon \Delta S(x)$ is between $S_0$ and $S_1$ (i.e. satisfying the constraint of the optimization).

$S_{opt}(x)$ is the optimal Seebeck profile implies

$$F\left[S_{opt}(x)+\varepsilon\Delta S(x)\right]\leq F\left[S_{opt}(x)\right] \text{ for } \varepsilon \text{ positive and small enough}$$

which leads to

$$\left.\frac{dF\left[S_{opt}(x)+\varepsilon\Delta S(x)\right]}{d\varepsilon}\right|_{\varepsilon=0} \leq 0$$

Notice that $\Delta S(x)$ satisfies $\int_{x_0}^{x_1}\Delta S(x)dx = 0$. We have

$$\left.\frac{dF\left[S_{opt}(x)+\varepsilon\Delta S(x)\right]}{d\varepsilon}\right|_{\varepsilon=0} = \frac{-2\left(\int_0^1 S_{opt}(x)dx\right)^2}{\left(\int_0^1 (1-x)S_{opt}^2(x)dx\right)^2} \times \int_0^1 (1-x)S_{opt}(x)\Delta S(x)dx$$

$$\Longrightarrow \quad \int_0^1 (1-x)S_{opt}(x)\Delta S(x)dx \geq 0$$

$$\Longrightarrow \quad \int_{x_0-\delta}^{x_0}(1-x)S_0\,dx - \int_{x_1}^{x_1+\delta}(1-x)S_1\,dx \geq 0$$

$$\Longrightarrow \quad \left(1-x_0+\frac{\delta}{2}\right)\delta S_0 - \left(1-x_1-\frac{\delta}{2}\right)\delta S_1 \geq 0$$

Dividing by $\delta$ and taking the limit as $\delta \to 0$, we obtain

$$(1-x_0)S_0 - (1-x_1)S_1 \geq 0$$

which contradicts with (A01). Therefore, $S_{opt}(x)$ must be continuous at both $x_0$ and $x_1$.

In other words, $S_{opt}(x)$ has the form

$$S_{opt}(x) = \begin{cases} S_0, & 0 \leq x \leq x_0 \\ \dfrac{q}{1-x}, & x_0 < x < x_1, \\ S_1, & x_1 \leq x \leq 1 \end{cases} \quad x_0 = 1-\frac{q}{S_0}, \quad x_1 = 1-\frac{q}{S_1} \quad \text{(A02)}$$

**Step #5: The optimal value of q is $q = S_0/2$**

$$\int_0^1 S_{opt}(x)dx = \int_0^{x_0} S_0 dx + \int_{x_0}^{x_1} \frac{q}{1-x}dx + \int_{x_1}^1 S_1 dx$$

$$= x_0 S_0 + q \cdot \log\left(\frac{1-x_0}{1-x_1}\right) + (1-x_1)S_1$$

$$= S_0 - q + q \cdot \log\frac{S_1}{S_0} + q$$

$$= S_0 + q \cdot \log\frac{S_1}{S_0}$$

$$\int_0^1 (1-x)S_{opt}^2(x)dx = \int_0^{x_0}(1-x)S_0^2 dx + \int_{x_0}^{x_1} \frac{q^2}{1-x}dx + \int_{x_1}^1 (1-x)S_1^2 dx$$

$$= \left(1 - \frac{x_0}{2}\right)x_0 S_0^2 + q^2 \cdot \log\left(\frac{1-x_0}{1-x_1}\right) + \frac{(1-x_1)^2}{2}S_1^2$$

$$= \frac{S_0^2 - q^2}{2} + q^2 \cdot \log\frac{S_1}{S_0} + \frac{q^2}{2}$$

$$= \frac{S_0^2}{2} + q^2 \cdot \log\frac{S_1}{S_0}$$

Consider the function

$$f(q) \equiv F[S_{opt}(x)] \equiv \frac{\frac{1}{2}\left(\int_0^1 S_{opt}(x)dx\right)^2}{\int_0^1 (1-x)S_{opt}^2(x)dx} = \frac{\left(S_0 + q \cdot \log\frac{S_1}{S_0}\right)^2}{S_0^2 + 2q^2 \cdot \log\frac{S_1}{S_0}}$$

Taking the derivative with respect to $q$, we have

$$\frac{df(q)}{dq} = \frac{2\left(S_0 + q \cdot \log\frac{S_1}{S_0}\right)\log\frac{S_1}{S_0}}{\left(S_0^2 + 2q^2 \cdot \log\frac{S_1}{S_0}\right)^2}\left[S_0^2 - 2S_0 q\right]$$

==> The maximum of $f(q)$ is attained at $q = \frac{S_0}{2}$.

In summary, the 5 steps above have completely determined the optimal solution

$$S_{opt}(x) = \begin{cases} S_0, & 0 \leq x \leq \dfrac{1}{2} \\ \dfrac{S_0}{2(1-x)}, & \dfrac{1}{2} < x < 1 - \dfrac{S_0}{2S_1} \\ S_1, & 1 - \dfrac{S_0}{2S_1} \leq x \leq 1 \end{cases}$$